\documentstyle[12pt]{article}
\advance\voffset by -3cm
\advance\hoffset by -3cm

\begin{document}
\title{Birational cobordisms and factorization of birational maps}


\author{\bf Jaros\l aw W\l odarczyk
\footnote{ The author is in part supported by Polish KBN Grant}
\date{Instytut Matematyki UW,\\  Banacha 2,
02-097 Warszawa, Poland,\\ jwlodar@mimuw.edu.pl }}

\maketitle

\bigskip





In this paper we develop a Morse-like theory in order to
decompose birational maps and morphisms of smooth
projective varieties defined over a field of characteristic
zero into more elementary
steps which are locally \'etale isomorphic to equivariant flips,
blow-ups and blow-downs of toric varieties (see Theorems 1,
2 and 3). The crucial role in the considerations is played by
$K^*$-actions where $K$ is the base field. The importance of
$K^*$-actions in birational geometry and their connection with Mori
Theory were already discovered by Thaddeus, Reid and many
others (see [Tha1], [Tha2], [Tha3], [R], [D,H]). On the
other hand, the ideas of the present paper were
inspired by the combinatorial techniques of  Morelli's proof of the
strong blow-up conjecture for toric varieties ([Mor]).

\bigskip

\noindent {\bf Acknowledgements}. I am  grateful to
Professor Andrzej Bia\l ynicki-Birula, Jaros\l aw Wi\'sniewski and
Jerzy Konarski from Warsaw University for numerous discussions and remarks.

\bigskip

The ground field $K$ is assumed to be algebraically closed.
 All algebraic varieties in
this paper and their morphisms are defined over $K$.

We recall the definitions of good and geometric quotients (see
also [Mum]).

\bigskip

\noindent {\bf Definition 1}. Let $K^*$ act on $X$. By a {\it good
quotient}  we mean a variety $Y=X//K^*$ together with
a morphism $\pi:X\rightarrow Y$ which is constant on
$G$-orbits such that for any affine open subset
$U\subset Y$ the inverse image $\pi^{-1}(U)$ is affine and
$\pi^*:O_Y(U)\rightarrow O_X(\pi^{-1}(U))^G$ is an isomorphism.
If additionally for any closed point $y\in Y$ its inverse limit $\pi^{-1}(x)$
is a single orbit we call $Y:=X/K^*$ together with
$\pi:X\rightarrow Y$ a {\it geometric quotient}.

\bigskip

\noindent {\bf Definition  2}: Let $X_1 $ and $X_2$ be two birationally equivalent
 normal varieties.
By a {\it birational cobordism} or simply a {\it cobordism}
$B:= B(X_1,X_2)$ between them we understand a
normal variety $B$ with an algebraic action of $K^*$ such that the
sets

\[\begin{array}{rccc}
B_-:&=\{x \in B\mid \, \lim_{t\to 0} \, tx \, \mbox{does
not exist}\}& \, \mbox{and}& \\
B_+:&=\{x \in B\mid\, \lim_{t\to\infty} \, tx \,
\mbox{does not exist}\} &&
\end{array}\]

\noindent are nonempty and open and there exist
 geometric quotients
$B_-/ K^*$ and $B_+\// K^*$ such that $B_+/ K^*\simeq X_1$ and $B_-\//
K^*\simeq X_2$ and the birational equivalence $X_1\mathrel{{-}\,{\rightarrow}}
X_2$ is given by the above
isomorphisms and the open embeddings
$ V:=B_+ \cap B_-\// K^*
\subset B_+\//K^*$ and $V\subset B_-\//K^*$.

\bigskip
\noindent {\bf Remark}. The analogous notion of cobordism of fans of
toric varieties was introduced by Morelli in [Mor].

\bigskip

 \noindent {\bf Remark}.
The above definition can  also be considered as an analog of the notion of
cobordism in Morse theory. In the present situation, however,
a Morse function defining in the classical theory
the local action of a 1-parameter group of diffeomorphisms, is
replaced by an action of $K^*$. The objects from Morse theory
like bottom and top boundaries, critical points can be interpreted
in terms of this action.

Let $W$ be a cobordism in Morse
theory of two
differentiable manifolds $X$ and $X'$ and $f:W\rightarrow
[a,b]\subset {\bf R}$ be a
Morse function such that $f^{-1}(a)=X$ and $f^{-1}(b)=X'$.
Then $X$ and $X'$ have open neighbourhoods $X\subseteq
V\subseteq W$ and $X'\subseteq
V'\subseteq W'$ such that  $V\simeq X\times
[a,a+\epsilon )$  and $V'\simeq X'\times (b-\epsilon ,b]$ for which
$f_{\mid V}:V\simeq X \times [a,a+\epsilon)\rightarrow [a,b]$ and
$f_{\mid V'}:V'\simeq X' \times (b-\epsilon ,b] \rightarrow [a,b]$ are
the natural projections on the second coordinate.  Let
$W':=W\cup_V X\times
(-\infty,a+\varepsilon) \cup_{V'}
X'\times (b-\varepsilon,+\infty )$. One can easily see
that $W'$ is isomorphic to $W\setminus X\setminus
X'=\{x\in W\mid  a < f(x) < b\}$. Let $f':W'\rightarrow {\bf R}$ be the
map defined by glueing the function $f$ and the natural projection on the second
coordinate. Then $grad (f')$ defines an action on $W'$ of
a $1$-parameter group
 $T\simeq {\bf R}\simeq {\bf
R}_{>0}^*$ of diffeomorphisms. The last group isomorphism is given by the exponential.

Set
\[\begin{array}{rc} W'_-:&=\{x \in W'\mid \, \lim_{t\to 0}
\,  tx  \, \mbox{does not exist}\},\\
W'_+:&=\{x
\in W\mid l\,im_{t\rightarrow\infty} \, tx \, \mbox{does
not exist}\}.
\end{array}\]
\noindent Then one can see that $W'_-$ and $W'_+$ are open and  $X$
 and $X'$ can be considered as quotients of these sets by $T$.
 The critical points of the Morse function are $T$-fixed points.

\bigskip

\noindent {\bf Example 1} (Atiyah [A]  and  Reid [R]).
Let $K^*$ act on the $(l+m)$-dimensional affine space
$B:=A^{l+m}_K$ by
$$ t(x_1,...,x_l,y_1,...,y_m)=(t\cdot x_1,...,t\cdot
x_l,t^{-1}\cdot y_1,...,t^{-1}\cdot y_m). $$

\noindent Set $\overline{x}:=(x_1,...,x_l), \, \overline{y}=(y_1,...,y_m) $.
Then
$$\displaylines { B_-= \{(\overline{x},\overline{y})\in A^{l+m}_K
\mid \, \overline{y}\neq 0\},\cr
B_+= \{(\overline{x},\overline{y})\in A^{l+m}_K \mid
\,\overline{x}\neq 0\}.\cr }$$

One can easily see that $B//K^*$ is the affine cone over
the Segre embedding of ${\bf P}^{l-1}\times {\bf P}^{m-1}\rightarrow
{\bf P}^{l\cdot m -1}$, and $B_+/K^*$ and $B_-/K^*$ are smooth.

The relevant birational map $\phi : B_-/K^* \mathrel{{-}\,{\rightarrow}} B_+/K^*$ is a
flip
for $l, m\geq 2$ replacing ${\bf P}^{l-1}\subset B_-/K^* $
with ${\bf P}^{m-1}\subset B_+/K^* $.
For $l=1, m\geq2$, $\phi$ is a blow-down, and for
$l\geq2, m=1$ it is a blow-up. If $l=m=1$ then  $\phi$ is the identity.

\bigskip

\noindent {\bf Remark}. In Morse theory we have an
analogous situation. In  cobordisms with one critical point we replace
$S^{l-1}$ by $S^{m-1}$.
\vskip6pt

The following example is a simple generalization of
Example 1.

\bigskip

\noindent {\bf Example 2} (Morelli  [Mor]).

Let $K^*$ act on $B:=A^{l+m+r}_K$ by
$$ t(x_1,...,x_l,y_1,...,y_m,z_1,...,z_r)=(t^{a_1}\cdot x_1,...,t^{a_l} \cdot
x_l,t^{-b_1}\cdot y_1,...,t^{-b_m}\cdot y_m, z_1,...,z_r).$$
where $a_1,...,a_l,b_1,...,b_m>0$.
Set $\overline{x}=(x_1,...,x_l) , \overline{y}=(y_1,...,y_m) ,
\, \overline{z}=( z_1,...,z_r)$.
Then
$$\displaylines{ B_-= \{(\overline{x},\overline{y},\overline{z})\in
A^{l+m+r}_K \mid  \, \overline{y}\neq 0\},\cr
B_+= \{(\overline{x}, \overline{y},\overline{z})\in
A^{l+m+r}_K \mid  \, \overline{x}\neq 0 \}\cr}$$

$B$, $B_-$ and $B_+$ can be considered as toric varieties
acted on by a torus
 $$T=\{(\overline{x},\overline{y},\overline{z})\in A^{l+m+r}_K
\mid \, x_i\neq 0  , y_j\neq 0  ,z_k\neq 0 \, \,  \mbox{for any} \, \,  i, j, k \}
.$$  This torus defines  the $ (l+m+r)$-dimensional lattices
$N:={\rm Hom}_{alg.gr}(K^*,T)$, $M:={\rm Hom}_{alg.gr}(T,K^*)$ and the
vector spaces $N_{\bf Q}=N \otimes {\bf
Q}$, $M_{\bf Q}=N\otimes {\bf Q}$. For any vectors $v\in N$ and
$w\in M$ let
$t_v$ denote the corresponding $1-$parameter subgroup and
$x_w$ denote the corresponding character. We have a perfect
pairing of lattices
$$\langle *,*\rangle:M\times N={\rm Hom}_{alg.gr}(T,K^*)\times
{\rm Hom}_{alg.gr}(K^*,T)\rightarrow {\rm Hom}_{alg.gr}(K^*,K^*)\simeq
{\bf Z}$$
given by $t^{\langle w,v\rangle}=x_w(t_v)$.

Then the action of $K^*$ on the relevant varieties determines
a $1$-parameter subgroup of $T$ which corresponds to a vector $v_0\in N$.
This defines a projection $\pi: N_{\bf Q}\rightarrow
N_{\bf Q}'$ where $N'_{\bf Q}:=N_{\bf Q}/({\bf Q}\cdot v_0)\supset N':=N/({\bf Q}\cdot v_0\cap N)$ . The dual vector space to $N_{\bf Q}'$
is  $M_{\bf Q}'=\{m\in M_{\bf Q}\mid  \langle m,v_0\rangle =0\}\supset
M':=\{m\in M \mid  \langle m,v_0\rangle =0\}$.

Consequently, $B\supset T$ is an affine toric variety corresponding to a regular
cone $\Delta\subset N_{\bf Q}$ and $B_-$ (respectively $B_+$) corresponds to the fan $\Delta_+$
(respectively $\Delta_-$) consisting of the
faces of $\Delta$ visible from above  (
respectively below).

Indeed,
$$\displaylines{p\in B_-\equiv\lim_{t\to 0}tp\,\,\mbox{does not
exist}  \cr \equiv \exists{ F \in
\Delta^{\vee}} \mbox{such that} \, \lim_{t\to 0}
x_F(tp)=\lim_{t\to 0} x_F(t_{v_0})x_F(p)\,\,\mbox{does not
exist} \,
\cr
\equiv \exists{  F \in
\Delta^{\vee}}\,\mbox{such that} \,
 \, \lim_{t\to
0}x_F(t_{v_0})=\lim_{t\to 0} t^{(F,{v_0})} \mbox{does not
exist} \,\mbox{and}
\cr
x_F(p)\neq 0 \,\equiv \exists{  F \in
\Delta^{\vee}} \, \,  \mbox{such that} \, \langle  F,{v_0}\rangle <0 \,
\mbox{and}\, p\in X_{\sigma_F }
\cr
\mbox{where} \, \, \sigma_F=\{v\in N_{\bf Q}\mid  \langle  F,v
\rangle =0\} \,\equiv\,
\, p \in
X_{\Delta_+}. \cr}$$
Analogously $B_+=X_{\Delta_-}$(see also [Jur], Thm. 1.5.3).

 The  quotients $B_+/K^*$ , $B_-/K^*$ and $B//K^*$
 are toric varieties corresponding  to the fans
$\pi(\Delta_+)=\{\pi(\sigma)\mid \sigma\in\Delta^+\}$,
$\pi(\Delta_-)=\{\pi(\sigma)\mid \sigma\in\Delta^-\}$ and $\pi(\Delta)$
respectively.

Indeed, $B//K^*=\mathop{\rm Spec} K[\Delta^{\vee}\cap M']$ corresponds
to the cone
$$(\Delta^{\vee}\cap M'_{\bf Q})^{\vee}=\{w\in
M'_{\bf Q}\mid \langle w ,* \rangle_{\mid (\Delta+{\bf Q}\cdot
v_0)}\geq 0\}^ {\vee}=
(\Delta+{\bf Q}\cdot v_0)/{\bf Q}\cdot v_0=\pi(\Delta)\subset
N'_{\bf Q}.$$
\noindent Analogously $B_-=X_{\Delta_+}$ (or $B_+=X_{\Delta_-}$)
is obtained by glueing together the affine pieces $X_{\pi(\sigma)}$ where
$\sigma\in \Delta_+$ (or respectively $\sigma\in
\Delta_-$). Obviously the projection $\Delta^+\rightarrow
\pi(\Delta^+)$ is 1-1 and the image
$\pi(\Delta^+)$ is a fan.

Note that $B_+/K^*$ and $B_-/K^*$ admit cyclic
singularities. This follows from the fact that $B_+/K^*$ is
covered by the open sets $$U_i:=\{(\overline{x}, \overline{y},\overline{z})\in
A^{l+m+r}_K \mid  \, x_i\neq 0 \}/K^* \simeq
\{(\overline{x}, \overline{y},\overline{z})\in
A^{l+m+r}_K \mid  \, x_i=1 \}/\Gamma$$ where $\Gamma =\{t\in
K^*\mid t(x_i)=t^{a_i}x_i = x_i\}$.

The relevant birational map $\phi : B_-/K^* \mathrel{{-}\,{\rightarrow}} B_+/K^*$
for $l,m\geq 2$ is a
toric flip associated with a bistellar operation replacing
the triangulation $\pi(\Delta_-)$ of the cone $\pi(\Delta)$
with  $\pi(\Delta_+)$ . It
 replaces the product of $A^r_K=\{ 0, 0,
\overline{z}) \in A^{l+m+r}_K\} $ and the $(l-1)$-dimensional weighted
projective space defined by the action of $T$ on the
$\overline{x}$-coordinates of $A^{l+m+r}_K$  with the product of $A^r_K$
and
the $(m-1)$-dimensional weighted
projective space defined by the action of $T$ on the
$\overline{y}$-coordinates  of $A^{l+m+r}_K$.
For $l=1, m\geq2$, $\phi$ is a toric blow-up whose
exceptional fibers are weighted projective spaces. For
$l\geq2,m=1$, $\phi$ is a toric blow-down. If $l=m=1$ then
$\phi$ is the identity.

\bigskip

\noindent {\bf Remark}. We prove in Theorems 1 and 2 that if $char K=0$ then any birational map of
any smooth projective or complete varieties can be decomposed into a sequence of
toroidal flips, toroidal blow-ups and blow-downs which
are locally \'etale isomorphic to toric flips, blow-ups, and
blow-downs described  in  Example 2.

\bigskip
\noindent {\bf Definition 3.}  Let $X_1 $ and $X_2$
be two birationally equivalent
 normal varieties and let $\varphi_1:X_1\rightarrow Y$
and $\varphi_2:X_2\rightarrow Y$ be two morphisms
commuting with the birational equivalence.
By a {\it birational cobordism over Y} \,
 between them we understand a
cobordism $B:= B(X_1,X_2)/Y)$ with  a $K^*$-equivariant morphism
$\phi :  B\rightarrow
Y$ where $Y$  is equipped
with the trivial $K^*$-action and such that the following
diagrams commute:

\[\begin{array}{rccc}
&B_+/K^* &\simeq &X_1 \\
&\uparrow&\searrow&\downarrow\\
\phi_{\mid B_+}:&B_+&\rightarrow&Y\\

\end{array}\]

\[\begin{array}{rccc}
&B_-/K^*&\simeq &X_2 \\
&\uparrow&\searrow&\downarrow\\
\phi_{\mid B_-}:&B_-&\rightarrow&Y\\
\end{array}\]

\noindent We say that the cobordism $B$ over $Y$ is
{\it trivial} over an open subset $U\subset Y$ iff
there exists an equivariant isomorphism $\phi^{-1}(U) \simeq
U\times K^*$, where the action of $K^*$ on $U\times K^*$ is
given by $t(x,s)=(x,ts)$.

\bigskip



We shall construct, in terms of the group action, an order
on the set of connected
components of the fixed point set which corresponds to the
order on the set of critical points defined by a Morse
function in Morse theory.

\bigskip

\noindent {\bf Definition 4}. Let $B$ be a cobordism. We say that
a connected component $F$ of the fixed point set  is {\it an
immediate predecessor} of a component $F'$ iff there exists a
non-fixed point $x$ such that $\lim_{t\to 0} tx\in F$ and
$\lim_{t\to \infty} tx\in F'$. We say that $F$  {\it precedes}
$F'$ and write $F<F'$ if there exists a sequence of connected
fixed point set components $F_0=F ,F_1,...,F_l=F'$ such that
$F_{i-1}$ precedes $F_i$ (see [B-B,S],  Def. 1.1). We call a
cobordism {\it collapsible} (see also [Mor]) iff the relation $<$
on its set of connected components of the fixed point set is an
order. (Here an order is just required to be transitive.)

\bigskip
\noindent {\bf Remark}. The concept of collapsibility in the
toric situation was introduced by Morelli in [Mor].

\bigskip

\noindent {\bf Defintion 5}. A cobordism $B$ is {\it projective} if
$B$ is a quasiprojective variety.


\bigskip

\noindent {\bf Lemma 1}. A projective  cobordism is
collapsible.

\bigskip

{\bf Proof}. By Sumihiro we can embed the variety
equivariantly into a projective space  ([Sum],  Thm. 1). Each connected
component of the fixed point set in the given variety
is contained in an irreducible component of the fixed point
set of the projective space. The homogeneous coordinates on
${\bf P}^n$ can be chosen to be semi-invariants in such a
way that

$$t([x_0,...,x_n])=[t^{a_0}x_0,...,t^{a_0}x_{l_0},t^{a_1}x_{l_0+1},...,
t^{a_k}x_{l_{k-1}+1},...,t^{a_k}x_{n}]$$,
\noindent where $0=a_0<a_1<...<a_k$.

The fixed point components have the following description:

$F_j=\{x\mid x_i=0 \,\, \mbox{for all}\, i \, \mbox{such that} \,
0 \leq i \leq l_{j-1} \, \mbox{and}
\, l_{j} <i \leq n\}$

for $j=0,...,k$. (We put here $l_1=-1$.)

We see that for any $a=[a_0,...,a_n]$,
$$\displaylines{\lim_{t\to 0} ta \in F_j$ iff $a_i = 0 \, \mbox{ for
all} \, i\leq l_{j-1} \mbox{and} \, \, \exists i\leq l_j \,\mbox{such
that}\, a_i \neq 0,\cr
\lim_{t\to \infty} ta \in F_j$ iff $a_i = 0 \, \mbox{ for
all} \, i> l_{j} \, \mbox{and} \, \, \exists i > l_{j-1} \, \mbox{such
that} \, a_i \neq 0 .\cr}$$

The order on the fixed point components on ${\bf P}^n$ is
determined by the relation:
$F_i<F_j$ iff $a_i<a_j$.

For any fixed point component $F$ on $X$ let $i(F)$ be the index
such that $F\subset F_{i(F)}$. The induced order on fixed point
components on $X$ satisfies: if $F<F'$ then $i_F<i_{F'}.$
The Lemma is proven.

\bigskip


\bigskip

Let $X$ be a variety with an action of $K^*$. Let
$F\subset X$ be a set of fixed points. Then we define  $$F^+(X)=F^+=\{x\in X\mid
\, \lim_{t\to 0} tx \in F\},\,  F^-(X)=F^-=\{x\in X\mid
\, \lim_{t\to \infty} tx \in F\}.$$

\bigskip

\noindent {\bf Definition 6}. Let $B$ be a collapsible
cobordism and $F_0$ be a minimal
component. By an {\it elementary collapse with respect to $F_0$}
we mean the cobordism $B^{F_0}:=B\setminus F_0^-$.
By an {\it elementary cobordism with respect to
$F_0$} we mean the cobordism $B_{F_0}:=B\setminus \bigcup_
{F\neq F_0} F^+$.

\bigskip

\noindent {\bf Proposition 1}. Let $F_0$ be a minimal component
 of the fixed point
set in a collapsible cobordism $B$. Then the
elementary collapse $B^{F_0}$ with respect to $F_0$ is again a
collapsible cobordism, in particular it satisfies:

a) $B^{F_0}_+=B_+$ is an open subset of $B$.

b) $F_0^-$ is a closed subset of $B$ and equivalently $B^{F_0}$
is an
 open subset of $B$.

c) $B^{F_0}_-$ is an open subset of $B^{F_0}$  and
$B^{F_0}_-=B^{F_0}\setminus \bigcup_{F\neq F_0} F^+$.

d) The elementary cobordism $B_{F_0}$ is  an open subset
of $B$ such that
$$\displaylines{{B_{F_0}}_-=B_{F_0} \setminus {F_0}^+=B_-\cr
{B_{F_0}}_+=B_{F_0} \setminus {F_0}^- =B^{F_0}_-.\cr}$$

e) There exist good and respectively geometric quotients  $B_{F_0}//{K^*}$
and $B^{F_0}_-/{K^*}$ and moreover the natural embeddings $ i_- : B_-\subset
B_{F_0}$ and $i_+ : B^{F_0}_-\subset B_{F_0}$ induce  proper morphisms
 $ i_{-/K^*} : B^{F_0}_-/{K^*}\rightarrow B_{F_0}//{K^*}$
 and  $ i_{+/K^*} : B_-/{K^*}\rightarrow B_{F_0}//{K^*}$.

\bigskip

Before proving the above proposition we shall state a few results
on $K^*$-actions.

\bigskip

\noindent {\bf Definition 7}. Let $X$ be a variety acted on
by
$K^*$. By a {\it sink} (resp. a {\it source}) of $X$ we mean an irreducible
component $F$ of the fixed point set such that $F^-$ (resp. $F^+$)
contains an open subset of $X$.

\bigskip

\noindent {\bf Lemma 2.} Let $X$ be a variety with a $K^*$-action for
which  $\lim_{t\to
\infty}tp$ exists for any $p\in X$. Then $X$ contains a sink $S$
and for any fixed point $x \in X$ there exists a sequence  $x_0=x,...,x_l, y_1,...,y_l$ of
points such that $\lim_{t\to
 0}ty_i=x_{i-1}$ and  $\lim_{t\to
\infty}ty_i=x_i$ for any $i=1,...,l$ and $x_l\in S$.

 \bigskip

\noindent {\bf Proof }. Let $X'$ be an equivariant
completion of $X$ (see [Sum], Thm.3). By ([Sum], Thm.2) we can find a projective
normal variety $X"$ with an action  of $K^*$ and an equivariant
birational morphism $\pi: X"\rightarrow X'$. Let $S"$ be a sink in $X"$.
Note that $\pi (S")$ is a sink in $X'$ and $\pi (S")\cap
X$  is a sink in $X'$ if it is non-empty. Let $x \in X$ be a
fixed point and let
$x' \in X"$ be a fixed point such that $\pi (x')=x$. By
([Sum], Thm 1)
we can embed $X"$ into a projective space and consequently
 find a sequence
$x'=x_0',...,x'_l, y'_1,...,y'_l$ in $X"$ as in the
statement of Lemma 2.
By applying the morphism $\pi$ we get a sequence
$x=x_0,...,x_l, y_1,...,y_l$
in $X'$. Note that if $x_i \in X$ then $y_{i+1} \in X$
since $X'$ is an open invariant neighbourhood of $x_i$ so
it contains all orbits having $x_i$  in their closure . On the
other hand by the assumption of Lemma 2 if $y_i \in X$ then
$x_{i+1}=\lim_{t\to\infty}ty_i \in X$. Finally the above sequence
is contained in $X$.

\bigskip

\noindent As a corollary from the proof of the above lemma we get

\noindent {\bf Lemma 3}.
Let $X$ be a variety with a $K^*$-action and with no sink.
 Then
 for any \\$y \in X$ there exists a sequence  $x_1,...,x_{l-1}, y=y_1,...,y_l$ of
points such that $\lim_{t\to
 0}ty_i=x_{i-1}$ for $i=2,...,l$, $\lim_{t\to
\infty}ty_i=x_i$ for any $i=1,...,{l-1}$ and $\lim_{t\to
 \infty}ty_l$ does not exist.

\bigskip

\noindent {\bf Proof of Proposition 1.}

a) This follows from the fact that $F_0^-\cap
B^+=\emptyset$.

b) Let $\overline{F_0^-}$ denote the closure of $F_0^-$ in
$B$. Let $\overline{F_0^-}=Z_0\cup...\cup Z_k$ be the
decomposition into irreducible components.
Since
$B_+$ is open and $F_0^- \cap B_+=\emptyset$ we see that
 $Z_i\cap
B_+=\emptyset$ for any $i$. In particular  $ \lim_{t\to
\infty}tp$ exists for any $p\in
Z_i$. It follows from  Lemma 2 that each $Z_i$ has a
sink $S_i$. Note that $S_i\subset F_0$. If not then by ([Kon],  Thm.
9),
$S_i^-(Z_i)$ contains an open subset $U_i\subset Z_i$ disjoint
from $F_0^-$. This gives
$\overline{F_0^-}\subset Z_0\cup...\cup (Z_i\setminus U_i)\cup...\cup Z_k$,
which is a contradiction.

On the other hand, if $Z_i\neq S_i^-(Z_i)$ then for any $x\in
Z_i\setminus S_i^-(Z_i)$ we can find a connected component $F$ of the fixed
point set of $B$ distinct from $F_0$ such that
  $\lim_{t\to
\infty}tx \in F$.
Now it follows from Lemma 2 applied to $Z_i$ that $F<F_0$,
which contradicts  the
assumption. Finally,
$\overline{F_0^-}\subset S_0^-(Z_0)\cup...\cup
S_k^-(Z_k)\subset F_0^-$, which means, that $F_0^-$ is closed.

c) We find directly from the definition of $B^{F_0}$ that
$B^{F_0}_- = B^{F_0}\setminus \bigcup_{F\neq F_0} F^+
= \ B \setminus (\bigcup_{F\neq F_0} F^+ \cup {F_0}^-$).

By repeating the argument in b) we see that
$\overline{ F^+} $  consists of points belonging to some
components $F'^+$ of the fixed point set such that $F'\geq F$. Hence
$B^{F_0}_-= B^{F_0}\setminus \bigcup_{F\neq F_0} \overline{ F^+}$ and
is open.

d) The same reasoning as above.

e) Set $U_0:= {B_{F_0}}\setminus F_0^-\setminus F_0^+ \subset
B_-$. Each orbit is closed in $U_0$. Since $B_-/K^*$
exists we deduce that $U_0/K^*$ exists.
Let $U_1,...,U_s$ be open affine invariant varieties
 covering $F_0$. Then $U_0,U_1,...,U_s$
cover $B_{F_0}$.

\noindent First we prove that

\vskip6pt
(*) each non-closed orbit is contained in some
closed invariant set of the form $\{x\}^+\cup \{x\}^-$ for
some $x \in F_0$.
\vskip6pt
\noindent This is equivalent to

$$\mbox{ for any}\,  y \in F_0, \,  \,
\{y\}^+\cup \{y\}^-\cap U_i \neq \emptyset \, \mbox{\, iff \,}
\, (\{y\}^+\cup
\{y\}^-)\subset U_i.$$  It suffices to prove that
$$ (U_i\cap F_0)^+ \cup
(U_i\cap F_0)^-=U_i\cap (F_0^+\cup F_0^-).$$

 Observe first that for any closed subset
$S\subset F_0$ the sets $S^-$ and $S^+$ are closed in $B_{F_0}$.
Let $\overline{S^-}$ be the closure of $S^-$ and let
$\overline{S^-}=Z_0\cup...\cup Z_k$ be its minimal decomposition
into irreducible components. As in the proof of b) we can show that
each $Z_i$ has  a sink $S_i$. Similarly if $Z_i\neq S_i^-(Z_i)$
then it follows from Lemma 2 that $F_0<F_0$, which is a
contradiction. Hence $Z_i=S_i^-(Z_i)$.
By the above we get
$$ \displaylines{S^-=(S_0\cap S)^-(Z_0)\cup...\cup (S_k\cap
S)^-(Z_k) \,  \mbox{and so} \cr \overline{S^-}=\overline{(S_0\cap
S)^-(Z_0)}\cup...\cup\overline{(S_k \cap S)^-(Z_k)}.\cr}$$ The last
equality implies $Z_i=\overline{(S_i \cap S)^-(Z_i)}$. This means
by ([Kon], Thm. 9) that $S_i\cap S= S_i$ or equivalently $S_i\subseteq
S$. Finally, $\overline{S^-}=(S_0)^-(Z_0)\cup...\cup
(S_k)^-(Z_k)\subseteq S^-,$ which means that $S^-$ is closed.

It follows  that $(U_i\cap F_0)^+\cup (U_i\cap F_0)^-
=F_0^+\cup F_0^- \setminus
(F_0\setminus U_i)^+ \setminus
(F_0\setminus U_i)^-$ is open in $(F_0^+\cup F_0^-)\cap U_i$.

Let $p_i:U_i\rightarrow U_i//K^*$ for $i=0,...,s$ denote the
standard projections. Then $p_i(U_i\cap F_0)$ are closed in
$U_i//K^*$, and thus $p_i^{-1}p_i(U_i\cap F_0)= (U_i\cap
F_0)^+\cup (U_i\cap F_0)^-$ are closed in $U_i$ which means
by the
connectedness of $F_0$ that
  $(U_i\cap F_0)^+\cup (U_i\cap F_0)^-=U_i\cap (F_0^+\cup F_0^-)$.
Thus (*) is proven.

Now since the fibers of $p_i$ are single orbits
which are closed in $B_{F_0}$ or are of the form $\{x\}^+\cup
\{x\}^-$ for $x\in F_0\cap U_i$ then for any $i\neq
j$ the map
$\phi_{ij}: U_i\cap U_j//K^*\rightarrow U_i//K^*$ is a
bijection and since the above varieties are normal it is an
open embedding. So we are in a position to define a good
quotient $B_{F_0}//K^*$ as a prevariety by glueing together
$U_i//K^*$ along
$(U_i\cap U_j)//K^*$. It suffices to prove that
$B_{F_0}//K^*$ is separated.

The closed embedding $F_0\subset B_{F_0}$ defines a map
$F_0\rightarrow B_{F_0}//K^*$. This map is bijective since each
fiber of the quotient map meeting $F_0$ is of the form
$\{x\}^+\cup \{x\}^-$ for $x \in F_0$. Moreover for any $i=0,1,...,s$
the map $F_0\cap U_i\rightarrow U_i//K^*$ is a closed embedding.
This implies that $F_0\rightarrow B_{F_0}//K^*$ is a closed embedding.
Let us identify $F_0$ with a
closed subset in $B_{F_0}//K^*$.

The open embedding $B_-\subset B_{F_0}$ defines a morphism
$\phi :B_-/K^*\rightarrow B_{F_0}//K^*$ whose restriction to
$(B_-/K^*)\setminus \phi^{-1}(F_0) $ is an isomorphism
$$\phi_{\mid (B_-/K^*)\setminus \phi^{-1}(F_0)} :(B_-/K^*)\setminus \phi^{-1}(F_0)\rightarrow
 (B_{F_0}//K^*)\setminus F_0.$$

Now let $R$ be any valuation ring and $K_0 \supset R$ be
its quotient field. Then we have the induced embedding
$\mathop{\rm Spec} K_0\hookrightarrow \mathop{\rm Spec} R$. In order to prove the separatedness of
$B_{F_0}//K^*$ we have to show that for any map
$f: \mathop{\rm Spec} K_0\rightarrow B_{F_0}//K^*$ there exists at most one
extension $f': \mathop{\rm Spec} R\rightarrow B_{F_0}//K^*$.

If $f(\mathop{\rm Spec} K_0)\subset F_0$ then $f'( \mathop{\rm Spec} R)\subset F_0$
and we are done by the separatedness of $F_0$.

If $f(\mathop{\rm Spec} K_0)\not\subset F_0$ then $f(\mathop{\rm Spec} K_0)\subset
(B_{F_0}//K^*)\setminus F_0\subset B_-/K^*$.

This gives the diagram
\[\begin{array}{rccccc}
& \mathop{\rm Spec} K_0&&\buildrel f \over  \longrightarrow&& B_-/K^*\\
&\downarrow&&\buildrel f_0\over  \nearrow&&\downarrow\\
 &\mathop{\rm Spec} R&& \buildrel f' \over  \longrightarrow&& B_F//K^*
\end{array}\]

It suffices to prove that each morphism $f'$ can be lifted
to a morphism $f_0$, in other words that the morphism $\phi$ is
proper. Then we are done by the separatedness of $B_-/K^*$.

The question is local so we can assume that $f'(\mathop{\rm Spec}R)\subset
U_i//K^* \subset B_{F_0}//K^*$ for some $i=0,...,s$. It
suffices to consider the case of $i\neq 0$ and $U_i$
affine.
 Let $\phi :U_i \hookrightarrow A^n_K$ be a closed embedding
into an affine space $A^n_K$ with a
linear action of $K^*$.
Let $F_A$ denote the fixed point set of $A^n_K$.

We have the following commutative diagram:
\[\begin{array}{rccccc}
&B_{F_0}//K^*& \supset &U_i//K^*& \hookrightarrow &A^n_K//K^*\\
&\uparrow  i_{-/K^*}& &\uparrow & &\uparrow t_- \\
&B_-/K^*& \supset & (U_i \setminus F_0^+) /K^*&
\hookrightarrow & (A^n_K\setminus
 F_A^+) /K^* \\
\end{array}\]
The morphism $t_-$ is proper and the horizontal arrows are closed
embeddings, hence $i_{-/K^*}$ is proper and $B_{F_0}//K^*$ is
separated.

Similarly since $B^{F_0}_-=B_{F_0+}$ is fixed point free we
see that $B_-/K^*$ is a prevariety. In order to prove the separatedness
of $B^{F_0}_-/K^*$ it is sufficient to prove the separatedness of the
morphism $B^{F_0}_-/K^*\rightarrow B_{F_0}//K^*$. We can reduce the
situation to the commutative diagram of separated
varieties
\[\begin{array}{rccccc}
&B_{F_0}//K^*& \supset &U_i//K^*& \hookrightarrow &A^n_K//K^*\\
&\uparrow i_{+/K^*}& &\uparrow & &\uparrow t_+\\
&B^{F_0}_-/K^*& \supset& (U_i \setminus F_0^-) /K^*&
\hookrightarrow & (A^n_K\setminus
F_A^-) /K^*\\
\end{array}\]

Properness of the relevant morphism follows from properness of
$t_+$.

\bigskip
As a corollary from Proposition 1 we get

\noindent {\bf Lemma 4}.  Let
$F_0,..., F_k$ be connected fixed point set components in a
collapsible cobordism $B$ such that
$F_i>F_j$ implies $i>j$ for any $i, j=1,...,k$. Then $B$ can be
represented as a union of elementary cobordisms
$$B=B_{F_0}\cup_{B_{F_0+}}B^{F_0}_{F_1}\cup_{B^{F_0}_{F_1+}}B^{F_0F_1
}_{F_2}\cup ...\cup B^{F_0...F_{k-2}}_{F_{k-1}}\cup_
{B^{F_0...F_{k-2}}_{F_{k-1}+}}B^{F_0...F_{k-1}}_{F_k} .$$

\bigskip

{\bf Construction of  a birational cobordism}.

\bigskip

Consider a line bundle $E$ over a normal variety $X$.
Let $s_E:X\rightarrow E$ be its zero
section. Then $E$ defines a line bundle $$E^\infty:=((E\setminus
s_E(X))\times ({\bf P}^1\setminus \{0\}))/K^*$$ where
the action of $K^*$ on $$(E\setminus
s_E(X))\times ({\bf P}^1\setminus \{0\})$$ is given by
$t(x,y)=(tx,t^{-1}y)$ for $x\in E\setminus
s_E(X)$ and $y\in {\bf P}^1\setminus \{0\}$. Here the action
of $K^*$ on $E$ is standard and the relevant action on
${\bf P}^1\setminus \{0\}$ is induced by the standard
embedding $K^*={\bf P}^1\setminus \{0\} \setminus
\{\infty\} \subset {\bf P}^1\setminus \{0\}).$

\bigskip

\noindent
{\bf Definition 8} ( [Nag2]).
 Let $X$ and $X'$ be
birationally equivalent  varieties with isomorphic open subsets
$X\supset U\simeq U'\subset X'$. Let $\Delta: U\rightarrow
X\times X'$ be the induced morphism. By the {\it join $X*X'$} of
$X$ and
 $X'$ we mean the closed subvariety
$ \overline{\Delta (U)}\subset X\times X'$ .

\bigskip

Now let $X\supseteq U\simeq U'\subseteq X'$ be  birationally
equivalent normal varieties. Let us identify
$U\simeq U'\simeq \Delta(U)$.
Denote by $\pi:X*X'\rightarrow X$ and
$\pi':X*X'\rightarrow X'$
the standard projections.

\bigskip

\noindent {\bf Lemma 5}. Let $D$, $D'$ be effective Cartier
divisors on $X$ and $X'$ respectively such that
$$V:=X*X'\,  \setminus \,  ( \pi^{-1}(supp(D))
\cup \pi'^{-1}( supp(D'))) \subseteq
 U. $$  Then  the open embeddings:

$$\displaylines{V\times K^*\subset O_X(-D)  \mbox{ and} \cr
V\times K^*\subset O_{X'}(D')^{\infty}. \cr}$$

\noindent(obtained by the natural multiplication by the sections
corresponding to $D$ and $D'$)
define
the separated set $$L(X,D;X',D'):=O_X(-D)\cup_{V\times K^*} O_X(D)^{\infty}.$$

\noindent {\bf Proof.}
Let $R$ be a valuation ring with  quotient field $K_0$ and
 valuation $\nu$. We have to prove that for any morphism
$\phi_0: \mathop{\rm Spec} K_0 \rightarrow L(X,D;X',D')$
we can find at most one morphism $\phi: \mathop{\rm Spec} R \rightarrow
L(X,D;X',D')$ which makes the following diagram commutative:
\[\begin{array}{rccccc}
&\mathop{\rm Spec} K_0&&\buildrel\phi_0 \over\longrightarrow& &L(X,D;X',D')\\
&\downarrow&&\buildrel\phi \over \nearrow&&\\
&\mathop{\rm Spec} R&&&&\\
\end{array}\]
If $\phi(\mathop{\rm Spec}K_0)\subset L(X,D;X',D')\setminus (V\times K^*)$ then we
are done since by  construction the last set is
obviously separated.

So we can assume that $\phi_0 (\mathop{\rm Spec} K_0) \subset V\times K^*$.
Let $V_0:=V,V_1,..., V_s$ (respectively $V'_0:=V,V'_1,...,
V'_{s'}$) be an open  covering of $ X$ (respectively of $ X'$)
such that for any $i:=1,...,l$ ($i':=1,...,l'$), $V_i$
(respectively $V'_{i'}$) is affine and  $D_{\mid V_i}$
(resp. $D'_{\mid V'_{i'}})$ is described by
$f_i\in K[V_i]$ (resp.$f'_{i'}\in K[V_{i'}]$).

Set $\psi_{ij}: (V_i\cap V_j)\times K\rightarrow V_i\times K,$
$\psi_{ij}(u,s_j)=(u,(f_j/f_i) s_j)$.
Here by $s_i$ we mean the standard coordinate function on $K=\mathop{\rm Spec}K[s_i]$.

Then $O_X(-D)$ is obtained by glueing $V_i\times K$ along
$\psi_{ij}$.

Analogously we can obtain $O_{X'}(D)^{\infty}$ by glueing together
$V'_{i'}\times ({\bf P^1}\setminus \{0\})$ along
$\psi'_{i'j'}$ where,
$\psi'_{i'j'}(u,s'_{j'})=(u,(f'_{j'}/f'_{i'})s'_{j'})$.

Here by $s'_{i'}$ we mean the standard coordinate function
on ${\bf P^1}\setminus \{0\}=\mathop{\rm Spec}K[1/s'_{i'}]$.

Let $t$ be the coordinate on $V\times K^*=V\times
\mathop{\rm Spec} K[t,t^{-1}]$.
Suppose we have two morphisms
$\phi,\phi': \mathop{\rm Spec} R\rightarrow
 L(X,D;X',D').$
Then we can assume that
$$\displaylines{\phi: \mathop{\rm Spec} R\rightarrow O_X(-D),\cr
\phi':\mathop{\rm Spec} R\rightarrow  O_{X'}(D')^{\infty}.\cr}$$

Assume that $\phi(\mathop{\rm Spec} R)\in
V_i\times K \subset O_X(-D)$ for some $i\neq 0$. Since
the functions $s_i=s_0/f_i=t/f_i$ and $f_i$ are regular on $V_i$ we have
$$\nu (\phi^*(t/f))=\nu (\phi_0^*(t/f))\geq 0.$$
Hence $\nu (\phi_0^*(t))\geq \nu (\phi_0^*(f))\geq 0.$

On the other hand, we can assume that $\phi'(\mathop{\rm Spec} R)\subset
V'_{i'}\times  ({\bf P^1}\setminus \{0\}) \subset O_{X'}(D')^{\infty}$ for some $i'\neq 0$. Then
since $1/s'_{i'}=1/(s'_0 f'_{i'})=1/(tf'_{i'})$ and $f'_{i'}$ are regular
on $V'_{i'}$ we have

$$\nu (\phi'^*(1/(f'_{i'}\cdot t)))=\nu (\phi_0^*(1/(f'_{i'}\cdot
t)))\geq 0.$$
\noindent Hence $\nu (\phi_0^*(t))\leq -\nu (\phi_0^*(f'_{i'}))\leq 0$.
Finally \, $ \nu ({\phi_0}^*(f_i))=\nu
({\phi_0}^*(f'_{i'}))= \nu ({\phi_0}^*(t))=0.$

Now let $p:O_X(-D)\rightarrow X$ and $p':O_{X'}(D')^{\infty}\rightarrow X'$
be the standard projections. Then $p \phi:\mathop{\rm Spec} R\rightarrow X$  and
$p' \phi':\mathop{\rm Spec} R\rightarrow X'$ define a morphism
$\overline{p \phi}:\mathop{\rm Spec} R\rightarrow X*X'$.

By the previous considerations and the assumptions
$$\overline{p \phi}(\mathop{\rm Spec} R) \subseteq
\{x\in V_i*V'_{i'}\mid f_i(x)\neq 0, f'_{i'}(x)\neq 0\}\subseteq V$$

\noindent and consequently $\phi (\mathop{\rm Spec} R)\subseteq V
\times K^*$  and $ \phi' (\mathop{\rm Spec} R) \subseteq V
\times K^*.$
But this contradicts the separatedness of $
V\times K^*$.

Lemma 5 is proven.

\bigskip

For a morphism $\phi: X\rightarrow Y$ and  a Cartier
divisor $D$ on $Y$ we denote by  $\phi^{*}(D)$  its inverse
transform. For any birational morphism $\phi: X\rightarrow Y$
of smooth varieties and a Weil divisor $D$ on $X$ we denote
by $\phi_{*}(D)$ its strict
transform.

\bigskip

\noindent
{\bf Lemma 6}.

\noindent {\bf A}. Let $X \supseteq U\simeq U'\subseteq X'$ be
normal and projective and $D$ and $D'$ be ample divisors on $X$ and $X'$
respectively such that
$$V:=X*X'\,  \setminus \,  ( \pi^{-1}(supp(D)) \cup \pi'^{-1}( supp(D'))) \subseteq
 U. $$
 Then $L(X,D;X',D')$ is
quasiprojective.

\noindent {\bf B}. Let $X \supseteq U\simeq U'\subseteq X'$ be
smooth and projective.
Assume that  $D$ is ample on $X$ and
$ supp(D) \supseteq X\setminus U$.
Then $L(X,D;X',0)$ is quasiprojective.

\noindent {\bf C}. Let $\phi : X\rightarrow X'$ be a  birational
morphism of smooth projective
varieties. Assume that $E$ is an effective divisor on $X$ such that
$-E$ is very ample relative to $X'$.
Let $D_X$ and $D_{X'}$ be very ample divisors  on $X$ and $X'$
respectively such that
 $D_X=\phi^*(D_{X'})-n\cdot E$ for some $n\in
{\bf N}$ (see [EGA], II, 4.6.13(ii)). Let $D:=D_X-\phi_*(D_{X'})=
\phi^*(D_{X'})-\phi_*(D_{X'}) - n\cdot E$.
Then $L(X,D;X',0)$ is quasiprojective.

\bigskip

{\bf Proof.} Let $ p: O_X(-D) \rightarrow X$ and $ p': O_{X'}(D')^{\infty}
\rightarrow X' $ be the natural projections.
Let $S_{0}\subset O_X(-D)$ be the zero section divisor
and
$S_{\infty}\subset O_{X'}(D')^{\infty}$ be
the infinity section divisor.

For any Weil divisor $D$ on an open subset $U\subset X$ let
$\overline{D}$ denote the closure of $D$ in $X$.

\bigskip

{\bf A}.
Set $D_0=p^*(D)$ and $D_1:=p'^*(D')$.

Then $D_0$ and $D_1$ are Cartier divisors on
$L(X,D;X',D')$ since their supports are closed in $L(X,D;X',D').$

Observe that

$S_{\infty}+{D_0} +(t) =
S_{0}+{D_1}. $

Find a natural number $n$ such that $nD$ and $nD'$ are very
ample divisors on $X$ and $X'$ respectively.
Now one can easily check that
$$D_L:=n{D_0}+nS_{\infty} \simeq n{D_1}+nS_0$$ is a base point free
divisor and  for any curve $C$ in $L(X,D;X',D')$ there exists an
 effective divisor equivalent to
$D_L$ which intersects
$C$. This means that $D_L$ defines a quasifinte morphism
$\phi:L(X,D;X',D')\rightarrow {\bf P}^n$,
which by the Zariski theorem ([Zar], [Mum2]) can be extended to a finite
morphism
$\overline{\phi}:\overline{L(X,D;X',D')}\rightarrow {\bf
P}^n$. Then
$\overline{\phi}^*(O(1))$ is ample on $\overline{L(X,D;X',D')}$ (see [Har2],
Prop. 4.4), which
means that $\overline{L(X,D;X',D')}$ is a
projective variety and $L(X,D;X',D')$ is quasiprojective.

\bigskip
{\bf B}.
Let $D'$ be ample  on $X'$.
Set $D_0=p^*(D)$ and $D_1=p'^*(D')$.
Then $S_{\infty}+\overline{D_0} +(t) =
S_{0} $. Again
find  $n$ such that $nD$ and $nD'$ are very
ample divisors on $X$ and $X'$ respectively.
Repeat the reasoning of case A for the divisor
$$D_L:=n{D_0}+n\overline{D_1}+nS_{\infty} \simeq n\overline{D_1}+nS_0.$$

\bigskip

{\bf C}.
Set $D_0:=p^*(D)$, $D_1=p^*(D_{X})$ and $D_2:=p'^*(D_{X'}).$
Then $\overline{D_1}=\overline{D_2}+D_0$ and $D_0+S_{\infty}+(t) = S_{0}.$

Analogously to case A we conclude that
$$D_L:=\overline{D_1}+S_{\infty}=\overline{D_2}+D_0+S_{\infty}
\simeq \overline{D_2}+S_{0}$$ is an ample divisor on the quasiprojective
variety $L(X,D;X',0)$.


\bigskip

{\bf Proposition 2}.
\noindent{\bf A}. There exists a projective
cobordism $B(X,X')$ between  any two  birationally equivalent normal
projective varieties $X$ and $X'$.

\noindent{\bf B}. There
exists a smooth
cobordism $B(X,X')$ between any   birationally equivalent smooth
varieties $X$ and $X'$ over a field of characteristic zero.

\noindent {\bf B'}. There
exists a smooth projective
cobordism $B(X,X')$ between any birationally equivalent
smooth projective
varieties $X$ and $X'$ over a field of characteristic zero.

\noindent{\bf C}. For any  birational morphism $X\rightarrow X'$  of smooth
projective varieties over a field of characteristic zero which
is an isomorphism  over $U\subset X'$, there
exists a smooth projective cobordism $B(X,X')/X'$ over $X'$  which
is trivial over $U$.

\bigskip

{\bf Proof.}

In cases A,
B', C one can find  divisors $D$ and $D'$ satisfying
respectively the conditions A, B, C of Lemma 6. In  case B we
find divisors $D$ and $D'$ satisfying the conditions of
Lemma 5.

\bigskip

{\bf A}. Let $\overline{L(X,D;X',D')}$ be
a $K^*$-equivariant projective completion of $L(X,D;X',D')$ (see
[Sum], Thm. 1). Let $\overline{B(X,D;X',D')}$ be its normalization.

{\bf B}. Let $\overline{L(X,D;X',0)}$ be
a $K^*$-equivariant  completion of  $L(X,D;X',0)$ (
[Sum], Thm. 3). Let $\overline{(B(X,D;X',0)}$
be its canonical  $K^*$-equivariant
resolution  (see [Hir] and [B-M]).

{\bf B'}. Let $\overline{L(X,D;X',0)}$ be
a $K^*$-equivariant projective completion of
$L(X,D;X',0)$. Let $\overline{B(X,D;X',0)}$ be its canonical $K^*$-equivariant
resolution .

{\bf C}. Let $\overline{L(X,D;X',0)}$ be
a $K^*$-equivariant projective completion of $L(X,D;X',0)$. Let
$L'$ denote the graph of the rational map
$\overline{L(X,D;X',D')}\rightarrow X'$ and
 $\overline{B(X,D;X',D')}$ be its canonical $K^*$-equivariant
resolution.

  Note that in all cases
$$L(X,D;X',D')\subset \overline{B(X,D;X',D')}.$$ (In cases
B, B' and C
we put $D'=0.$)

Let $S_{0}\subset O_X(-D)\subset L(X,D;X',D')$ be the zero section divisor
and
$S_{\infty}\subset O_{X'}(D')^{\infty}\subset L(X,D;X',D')$ be
the infinity section divisor.

Set $$B(X,X'):=\overline{B(X,D;X',D')}\setminus S_{0} \setminus S_{\infty}.$$
\noindent Then $B(X,X')_+= \{x\in B(X,X')\mid \lim_{t\to\infty} \, tx
\, {\rm does\,
not \, exist}\}  =
\{x\in B(X,X')\subset \overline{B(X,D;X',D')}:
 \lim_{t\to\infty}tx\in
S_{0}\cup S_{\infty}\} = \{x\in B(X,X')\mid  \lim_{t\to\infty}tx\in
S_{\infty}\} = B(X,X')\cap S_{\infty}^- =
O_{X'}(D')^{\infty}\setminus S_{\infty}.$

Analogously $B_-(X,X')=O_X(-D)\setminus S_{0}$

In both cases evidently $B_+/K^*\simeq X'$ and
$B_-/K^*\simeq X$.

\bigskip

\noindent {\bf Remark}. The above constructed
cobordism $B$ between $X$ and $X'$ is of the form
$\overline{B}\setminus X'\setminus X$ where $\overline{B}$ is a complete variety
with a $K^*$-action such that $X$ is its source and $X'$ is
its sink. One can prove that each cobordism is of that form
(see Lemma 7).
This makes the analogy between birational cobordism and
cobordism in Morse theory stronger.

\bigskip

Another  method of constructing cobordisms in case C
was communicated to me by Abramovich.

Let $X_2\rightarrow X_1$ be a projective morphism of two smooth
varieties. Let $I\subset O$ be a sheaf of ideals such that
 $X_2=Bl_I{X_1}$  is obtained from $X_1$ by blowing up of $I$.
 Let
$W=X_1\times {\bf P}^1$ and let $\pi_1:W\rightarrow X_1$,
$\pi_2:W\rightarrow {\bf P}^1$ be the standard projections.
Let $z$ denote the standard coordinate on ${\bf P}^1$ and
let $I_0$ be the ideal of the point $z=0$ on ${\bf P}^1$. Then
$I'={\pi_1}^*(I)+{\pi_2}^*(I_0)$ is an ideal supported on $X_1\times \{0\}$.
Set $W'=Bl_{I'}W$.
 Then the proper transform
$\overline{X_1\times\{0\}}$ of $X_1\times
\{0\}$ is isomorphic to $X_2$ and it is a sink of $W'$. For
simplicity we identify it with $X_2$.  We prove
that $W'$ is smooth at $X_2\subset W'$.  Let $f_1,...,f_k,
z$ generate the ideal $I'$ at $x\in X_1\times
\{0\}.$
Then the completion of local ring of any point $y$ of $X_2$ is up
to linear transform of $f_1,...,f_k$
 equal to $ \widehat{O_y}=\widehat{O_x}[[f_1,f_2/f_1,..., f_k/f_1, z/f_1]]$ where $z/f_1$
generates the ideal of $X_2$. Since we know that $X_2$ is
smooth we find that
$\widehat{O_x}[[f_1,f_2/f_1,..., f_k/f_1]]$ is regular and finally
since $ z/f_1$ is algebraically independent of elements of
$\widehat{O_x}[[f_1,f_2/f_1,..., f_k/f_1]] $  we conclude that that $\widehat{O_y}$ is
regular which gives the smoothness of $W'$ at $X_2$. Now it is
sufficient to apply the canonical resolution and we get a smooth
variety with  sink $X_2$ and  source $X_1$. By [B-B] we conclude that
$X_2^+$ is a locally trivial $K$-bundle and finally
$W'\setminus X_1 \setminus X_2$
is a smooth projective cobordism from $X_2$
to $X_1$.

 \bigskip

\noindent {\bf Lemma 7}.
{\bf A}. Let $B$ be a normal variety with a $K^*$-action
with no fixed
points and such that the geometric quotient   $B/K^*$
exists. Then there exists a normal variety \\$B^0=B\,\cup \,
(B/K^*)$ (respectively $B^\infty = B\,\cup \, (B/K^*)$) with
a $K^*$-action such
that $B^0//K^*\simeq B/K^* \subset {B}^0$ is a source
in $B^0$ (respectively $B^\infty//K^*\simeq B/K^* \subset
B^\infty$ is a sink in  $B^\infty$ ) and the standard projection
$B^0\rightarrow B^0//K^*$ (resp. $B^\infty\rightarrow B^\infty//K^*$)
is given by  $x\in B^0\longrightarrow \lim_{t\to 0}
tx$  ( $x\in B^\infty\longrightarrow \lim_{t\to
\infty} tx$).

{\bf B}. Let $B(X,X')$ be a cobordism between $X$ and $X'$. Then
there exists a variety $\overline{B(X,X')}=B(X,X')\cup X\cup X'=
B(X,X')\cup_{B(X,X')_+}(B(X,X')_+)^\infty\cup_{B(X,X')_-} (B(X,X')_-)^0
$
with  source $X$ and  a sink $X'$.
If $X$ or $X'$ is complete than $\overline{B(X,X')}$ is also
complete.

\bigskip

\noindent {\bf Proof}. {\bf A}.  Let $K^*$ act on
$B\times K$ by $t(x,s):=(tx,t^{-1}s)$ where $x\in B$
and $s\in K$.
 This action is fixed point free and consequently
the quotient
$B^0:=(B\times K)/K^*$ is a prevariety. The morphism
$(B\times K)/K^*\rightarrow B/K^*$ is
separated since its restriction to any open affine invariant
$U\subset B$ determines a separated morphism $(U\times
K)/K^*\rightarrow U/K^*$. This implies the separatedness of
$B^0$.
The quotient $(B\times K^*)/K^*$
is isomorphic to $B$.  The morphism $i:B\simeq (B\times
K^*)/K^* \rightarrow (B\times
K)/K^*=B^0 $ is a 1-1 morphism of normal varieties and hence it is an
open embedding. The action of $K^*$ on $B\times K $ defined by
$t(x,s)=(tx,s)$ or equivalently $t(x,s)=(x,ts)$ induces an
action on $B^0$ which extends the given action on $B$.

Moreover
$$(B \times K)/K^* \/   \setminus \/ i((B\times
K^*)/K^*)= B^0\setminus B =(B \times {0})//K^* \simeq B//K^*. $$

Let $\pi: B\times K\rightarrow (B\times K)/K^*$ be the
standard projection. Let $\overline{x}=\pi(x,s) \in
(B\times K)/K^*$. Then $\lim_{t\to 0}
t\overline{x}=\pi(\lim_{t\to 0}(x,st))=\pi(x,0)\in (B/K^*)\times \{0\}$.

  The above reasoning can be repeated for
$B^\infty:=(B\times({\bf P}^1\setminus\{0\}))/K^*$.

\bigskip

\noindent {\bf B}. By A we can construct the prevariety
$\overline{B(X,X')}=(B_+)^\infty\cup_{B_+} B \cup_{B_-} (B_-)^0$ as in the
statement of the Lemma. We prove that this prevariety is
separated. We show first that $(B_+)^\infty\cup_{B_+}B$ is
separated. It is sufficient to show
that $\Delta(B_+)\subset (B_+)^\infty\times B$ is closed.
Let $\pi_1:\overline{\Delta(B_+)}\rightarrow (B_+)^\infty$
and $\pi_2:\overline{\Delta(B_+)}\rightarrow B$ denote the standard
projections. They are birational morphisms which are
isomorphisms over $B_+$.
Suppose on the contrary that there exists $y\in
\overline{\Delta(B_+)}\setminus \Delta(B_+)$. In particular
$\pi_1(y)\in (B_+)^\infty\setminus B_+)$ belongs to the
sink in $(B_+)^\infty$.
The points $p$ from the set ${\Delta(B_+)}\subset
\overline{\Delta(B_+)}$ have neither of the two limits $\lim_{t\to
0} tp$ and $\lim_{t\to \infty} tp$ in
$\overline{\Delta(B_+)}$. In particular $\overline{\Delta(B_+)}$
has no sink. By Lemma 3 we find a sequence
 $x_1,...,x_{l-1}, y=y_1,...,y_l$ of points in $\overline{\Delta(B_+)}$
such that $\lim_{t\to
 0}ty_i=x_{i-1}$ for $i=2,...,l$ and $\lim_{t\to
\infty} ty_i=x_i$ for  $i=1,...,l-1$ and $\lim_{t\to
\infty}ty_l$ does not exist.
By the previous remark $\pi_1(y_l)$  also belongs to a sink of
$(B_+)^\infty$. This means that $\lim_{t\to
\infty}t\pi_2(y_l)$ does not exist. But this
implies that $\pi_2(y_l)\in B_+$ and finally $y_l\in
\Delta(B_+)$ and $\pi_1(y_l)\in B_+$, a contradiction. We
have proved
that $(B_+)^\infty\cup_{B_+}B$ is separated.

Similarly one can prove that $\overline{B(X,X')}=((B_+)^\infty\cup_{B_+} B)
\cup_{B_-} (B_-)^0$ is separated.

Note that
$\lim_{t\to 0} tx$ and $\lim_{t\to \infty} tx$
exist for any $x \in \overline{B(X,X')}$. Now assume that $X'$ (or $X$) is complete. Let $\overline{B(X,X')}'$
be the normalization of a completion of $\overline{B(X,X')}$ . Then
$X$ is a sink in $\overline{B(X,X')}'$. By Lemma 2 for any
$x\in\overline{B(X,X')}'$ we can find a sequence
$x_0=x,...,x_l, y_1,...,y_l$ in  $\overline{B(X,X')}'$ such that $\lim_{t\to
 0}ty_i=x_{i-1}$,  $\lim_{t\to
\infty}ty_i=x_i$ for any $i=1,...,l$ and $x_l\in X$. As in the
proof of Lemma 2 we can show that $x\in  \overline{B(X,X')}$, which means
that $\overline{B(X,X')}=\overline{B(X,X')}'$ is complete.

\bigskip

\noindent {\bf Lemma 8}. Let $B_{F_0}$ be a smooth elementary cobordism. Then for
any $x\in F_0$ there exists an invariant neighbourhood
$V_x$ of $x$ and a
$K^*$-equivariant \'etale morphism $\phi :V_x\rightarrow
T_x $,
where $T_x\simeq A_k^n$ is the tangent space with the
induced linear $K^*$-action, such that in the diagram
\[\begin{array}{rccccc}
&V_x//K^* \times_{ T_x//K^*}
{T_x}_-/K^*& \simeq &{V_x}_-/K^*& \rightarrow
&{T_x}_-/K^* \\

&&& \downarrow & &\downarrow \\

&&&V_x//K^*    & \rightarrow  &T_x//K^* \\
&&&\uparrow& &\uparrow \\
&V_x//K^*  \times_  { T_x//K^*}
{T_x}_+/K^* & \simeq &{V_x}_+/K^*&\rightarrow & {T_x}_+/K^*\\

\end{array}\]

\noindent the vertical arrows are defined by open embeddings
and the horizontal morphisms are defined by $\phi$ and are \'etale.

\bigskip

\noindent {\bf Proof}.
By taking local semi-invariant parameters at the point $x\in
F_0$ one
can construct an equivariant morphism $\phi:U_x\rightarrow
T_x\simeq A^n_K$ from some
open affine invariant neighbourhood $U_x$ such that $\phi$ is \'etale
at $x$.
By Luna's Lemma (see [Lu], Lemme 3 (Lemme Fondamental))
there exists an invariant affine
neighbourhood $V_x\subseteq U_x$ of the point $x$ such that
$\phi_{\mid V_x}$ is \'etale, the induced map
$\phi_{\mid V_x/K^*}:V_x/K^* \rightarrow T_x/K^*$ is \'etale
and $V_x\simeq V_x/K^* \times_{T_x/K^*} T_x$.
It follows from the last property that $\phi_{\mid Gy}$ is an
embedding for any $y\in V_{x+}$.

Now
since $V_{x+}$ is an open invariant subset of the fixed point
free $B_+$ and $B_+/K^*$ exists it follows that $V_{x+}/K^*$ also exists.
Again by the Luna Lemma applied to affine neighbourhoods of
any $y\in V_{x+}$ and $\phi(y) \in T_x $ we deduce that
$ \phi_{/K\mid V_{x+}/K^*}:V_{x+}/K^*
\rightarrow T_x /K^*$ is \'etale at each point $y\in
V_{x+}/K^*$ and consequently is \'etale.

Let $\psi_+:{V_x}_+/K^* \rightarrow V_x//K^*  \times_  { T_x//K^*}
{T_x}_+/K^*$ be the natural map. Let $\pi_1: V_x//K^*  \times_  { T_x//K^*}
{T_x}_+/K^*\rightarrow {T_x}_+/K^*$ and $\pi_2: V_x//K^*  \times_  { T_x//K^*}
{T_x}_+/K^*\rightarrow V_x//K^* $ denote the natural
projections. Since $ V_x//K^* \rightarrow T_x//K^*$ is
\'etale we infer that $\pi_1$ is \'etale. On the other hand,
$\phi_{/K\mid V_{x+}/K^*}= \pi_1 \psi_+$
is \'etale. This implies that $\psi_+$ is \'etale and in particular
quasifinite. Since it is
a birational quasifinite morphism of normal varieties it is an
open embedding by the Zariski Main Theorem ([Mum2]).

On the other hand, $i_{/K^*}: {V_x}_+/K^* \rightarrow
V_x//K^*$ is proper by Proposition 1. But $i_{/K^*}=\pi_2
\psi_+$ and since $\pi_+$ is separated we conclude that $\psi_+$
is proper (see [Har1], Cor. 4.8 e)). Finally, a proper morphism which is an
open embedding is an isomorphism.

\bigskip

\noindent {\bf Definition 9} (see also [Dan]). A variety $X$ is called {\it
toroidal} iff for any $p\in X$ there exists an open affine
neighbourhood $U_x$ and an \'etale map $\phi:U_x\rightarrow
X_{\sigma_x}$ into some affine toric variety
$X_{\sigma_x}$. $X$ is called {\it quasismooth toroidal} iff
${\sigma_x}$ is a simplicial cone for any $x\in X$.

\bigskip
\noindent {\bf Definition 10} (see also [Oda], [Mor], [Wlo]).
 Let $\sigma=\langle  v_1,...,v_r\rangle \,
\subseteq N_{\bf Q}:= N\otimes {\bf Q}\simeq {\bf Q}^k$ be
an r-dimensional simplicial cone spanned by integral
vectors $v_1,...,v_r \in N\simeq {\bf Z}^k  $. Let
$\tau=\langle v_1,...,v_s\rangle\, \subset \sigma$ be
its face and let $\rho \in {\rm Relint}(\tau)$ and let $v_\rho$ be the
generator of $N\cap \rho$. By a {\it star subdivision}
$ \sigma_\rho$ of $\sigma$ at $\rho$ we mean the fan whose
set of maximal cones is  \\ $$\{\langle v_\rho,v_1,...,v_{i-1},
\check{v}_i, v_{i+1},...,v_s,...,v_r\rangle
\mid  \, 1\leq i \leq s\}.$$
We call the corresponding birational toric morphism
$X_{\sigma_\rho}\rightarrow X_\sigma $ a {\it toric blow-up} of $X_\sigma$.

\bigskip
\noindent {\bf Definition 11}. Let $X$ and $Y$ be two
quasismooth toroidal varieties. Then a birational morphism
$\psi : X\rightarrow Y$ is called {\it a toroidal blow-up} iff
the basic set $L\subset Y$ of $\psi$ is irreducible
 and there
is a simplicial cone $\tau \subset N_{\bf Q}$ and a ray $\rho \in
{\rm Relint} (\tau)$
such that for any $y\in L$  there
exists an open neighbourhood $U_y$ and a  commutative diagram

\[\begin{array}{rccccc}
&& & U_y & \rightarrow & X_{\sigma_y}\\

&&&\uparrow \psi & & \uparrow bl_{X_{\sigma_y}} \\

U_y \times_{X_{\sigma_y}}
X_{\sigma_{y,\rho}} &&\simeq&{\psi^{-1}(U_y)}
 & \rightarrow &
X_{\sigma_{y,\rho}} \\

\end{array}\]

\noindent where all the horizontal arrows are \'etale,
   $\tau$ is a face of $\sigma_y \subset N_{\bf Q}$ and
$bl_{X_{\sigma_y}}: X_{\sigma_{y,\rho}} \rightarrow
X_{\sigma_{y}}$ is a toric blow-up.

If $\sigma_y=\tau$ for any $y\in L$ then $\phi_y$ is called a {\it simple
toroidal blow-up}.

\bigskip

\noindent {\bf Remark}. It follows from the definition that
the basic set of a toroidal blow-up is a quasismooth
toroidal variety, and the basic set of a simple toroidal
blow-up is  smooth. The exceptional divisor of a simple
toroidal blow-up is a locally
free bundle whose fibers are toric varieties associated
with the fan $\{\pi(\tau')\mid  \, \tau \, \mbox{is a proper face of}
\, \tau \}$ where $\pi:N_{\bf Q}\rightarrow N_{\bf Q}/{\bf Q}\rho$ denotes the
standard projection. In particular if $Y$ is smooth
the fibers are just weighted projective spaces.

\bigskip
\noindent {\bf Defintion 12} (see also [Mor], [Wlo]). Let
$\sigma=\langle v_1,...,v_{k+1}\rangle$ be a
$k$-dimensional cone generated by $k+1$ integral vectors
$v_1,...,v_{k+1}$ with a unique relation $\sum a_iv_i=0$
where $a_i>0$ for all $1\leq i\leq l$ and $a_i<0$ for
$l+1\leq i\leq k+1$ where $l$ is some number $2\leq l\leq k$.
By a {\it stellar transform} of $\sigma$ we mean the
transformation replacing the subdivision $\Sigma_1$ with the set of
maximal simplices $\{\langle  v_1,...,\check{v_i},...,v_{k+1}\rangle \mid
1\leq i\leq l\}$ with  another subdivision  $\Sigma_2$ with the set of
maximal simplices $\{\langle v_1,...,\check{v_i},...,v_{k+1}\rangle \mid
l+1\leq i\leq k+1\}$.

We call the  diagram
\[\begin{array}{rccccc}
&X_{\Sigma_1}&&&&X_{\Sigma_2}\\
&&\searrow&& \swarrow&\\
&&&X_{\sigma}&&\\
\end{array}\]

\noindent a {\it toric flip}.

\bigskip

\noindent {\bf Definition 13}. Let $X$, $Y$ be quasismooth toroidal
varieties and $Z$ be any toroidal variety. Then we call a commutative diagram
 \[\begin{array}{rccccc}
&X&&&&Y\\
&&\psi_X\searrow&& \swarrow\psi_Y&\\
&&&Z&&\\
\end{array}\]

\noindent a {\it simple toroidal flip} if $X$ and $Y$ are quasismooth
toroidal varieties and $Z$ is a toroidal variety such that
the basic sets $L_X\subset Z$ and $L_Y\subset Z$ of $\psi_X$ and
$\psi_Y$ respectively  coincide and are irreducible
and
there exists a toric flip
\[\begin{array}{rccccc}
&X_{\Sigma_1}&&&&X_{\Sigma_2}\\
&&\searrow&& \swarrow&\\
&&&X_{\sigma}&&\\
\end{array}\]
such that
for any $z$ in the basic set $L_X=L_Y$ there
exists an open neighbourhood $U_z\subset Z$ and a commutative diagram
of morphisms
\[\begin{array}{rccc}

U_z \times_{X_{\sigma}} X_{\Sigma_1} \simeq
&{\psi_X^{-1}(U_z)}&\rightarrow &
X_{\Sigma_1} \\
&\downarrow \psi_X & & \downarrow \\
 & U_z & \rightarrow & X_{\sigma}\\
&\uparrow \psi_Y & & \uparrow  \\
U_z \times_{X_{\sigma}} X_{\Sigma_2} \simeq
&{\psi_Y^{-1}(U_z)}& \rightarrow &
X_{\Sigma_2} \\

\end{array}\]
where all the horizontal arrows are \'etale

\bigskip

\noindent {\bf Remark}. It follows from the definition that
the basic set $L_X=L_Y$ of  $\psi_X$ and $\psi_Y$ is smooth.
The exceptional divisors of $\psi_X$ and $\psi_Y$ are locally
free bundles whose fibers are quasismooth toric varieties.
In particular if $X$ and $Y$ are smooth  the
fibers are  weighted projective spaces .

\bigskip

\noindent {\bf Lemma 9}. Let $B_{F_0}$ be a smooth elementary cobordism.
Then the diagram  \[\begin{array}{rccccc}
&B_{F_0-}/K^*&&&&B_{F_0+}/K^*\\
&&\psi_-\searrow&& \swarrow\psi_+&\\
&&&B_{F_0}//K^*&&\\
\end{array}\]
is either
a simple toroidal flip such that the fibers of $\psi_-$ and
of $\psi_+$ are weighted projective spaces
or

$\bullet\psi_-$ is an isomorphism and $\psi_+$ is a simple toroidal blow-up
whose fibers are  weighted projective spaces
or

$\bullet\psi_+$ is an isomorphism and $\psi_-$ is a simple toroidal blow-up
whose fibers are  weighted projective spaces.

\noindent {\bf Proof}. For any $x\in F_0\subset B_{F_0}$ the
semi-invariant local parameters at $x$ determine a linear action on
tangent space at $x$. These tangent spaces  are equivariantly
isomorphic for all points of $F_0$ and determine a unique (up
to isomorphism)  affine
space with a linear action. It is sufficient to apply Lemma 8
and cite Example 2.

\bigskip

\noindent {\bf Lemma 10}. The birational equivalence
determined by a simple toroidal flip
\[\begin{array}{rccccc}
&X&&&&Y\\
&&\psi_X\searrow&& \swarrow\psi_Y&\\
&&&Z&&\\
\end{array}\]

is the
composite of a toroidal blow-up $\widetilde{X\times_Z Y}\rightarrow X$
and a toroidal blow-down  $Y\leftarrow \widetilde{X\times_Z
Y}$, where $\widetilde{X\times_Z Y}$ is the
normalization of $X\times_Z Y$

\noindent {\bf Proof}. Take $z\in L_X=L_Y$. The
commutative diagram

\[\begin{array}{rccccc}
&&&\widetilde{X\times_Z Y}&&\\
&&&\downarrow&&\\
&&&X\times_Z Y&&\\
&&\swarrow&& \searrow&\\
&X&&&&Y\\
&&\searrow&& \swarrow&\\
&&&Z&&\\
\end{array}\]

defines locally a diagram
\[\begin{array}{rccccccc}
&&&\widetilde{{\psi_X^{-1}(U_z)}\times_Z {\psi_Y^{-1}(U_z)}}&&&&\\
&&&\downarrow&&&&\\
&&&{\psi_X^{-1}(U_z)}\times_Z {\psi_Y^{-1}(U_z)}&&&&\\
&&\swarrow&& \searrow&&&\\
&\psi_X^{-1}(U_z)&&&&\psi_Y^{-1}(U_z)&&(*)\\
&&\searrow&& \swarrow&&&  \\
&&&U_z&&&&\\
\end{array}\]

On the other hand consider the  diagram of toric varieties
\[\begin{array}{rccccc}
&&&\widetilde{X_{\Sigma_1}\times_{X_{\sigma}}X_{\Sigma_2}}&&\\
&&&\downarrow&&\\
&&&X_{\Sigma_1}\times_{X_{\sigma}}X_{\Sigma_2}&&\\
&&\swarrow&& \searrow&\\
&X_{\Sigma_1}&&&&X_{\Sigma_2}\\
 &&\searrow&& \swarrow&\\
&&&X_{\sigma}&&\\
\end{array}\]

\noindent It follows from the universal property of the  fiber product that
$\widetilde{X_{\Sigma_1}\times_{X_{\sigma}}X_{\Sigma_2}}$ is a
normal toric variety whose fan consists of the cones $
\{\tau_1\cap\tau_2 \mid  \tau_1\in \Sigma_1, \tau_2\in \Sigma_2\}.$
The morphisms
$\widetilde{X_{\Sigma_1}\times_{X_{\sigma}}X_{\Sigma_2}}
\rightarrow X_{\Sigma_1}$ and $\widetilde{X_{\Sigma_1}\times_{X_{\sigma}}X_{\Sigma_2}}
\rightarrow X_{\Sigma_2}$  are toric blow-ups.
Moreover the above diagram induces the following one:

\[\begin{array}{rccccc}
&&&\widetilde{X_{\Sigma_1}\times_{X_{\sigma}}X_{\Sigma_2}}
\times_{X_{\Sigma_1}\times_{X_{\sigma}}X_{\Sigma_2}}{\psi_X^{-1}(U_z)}
\times_Z {\psi_Y^{-1}(U_z)}&&\\
&&&\downarrow&&\\
&&&{\psi_X^{-1}(U_z)}\times_Z {\psi_Y^{-1}(U_z)}&&\\
&&\swarrow&& \searrow&\\
&\psi_X^{-1}(U_z)&&&&\psi_Y^{-1}(U_z)\\
&&\searrow&& \swarrow&  \\
&&&U_z&&\\
\end{array}\]

Now it is sufficient to show that the above diagram coincides
with $(*)$. To this end we note that the morphism
$$\widetilde{{\psi_X^{-1}(U_z)}\times_Z
{\psi_Y^{-1}(U_z)}}\rightarrow \widetilde{X_{\Sigma_1}\times_{X_{\sigma}}X_{\Sigma_2}}
\times_{X_{\Sigma_1}\times_{X_{\sigma}}X_{\Sigma_2}}{\psi_X^{-1}(U_z)}
\times_Z {\psi_Y^{-1}(U_z)}$$ is proper birational and \'etale.
Both varieties are normal since the completions of local rings
of the second variety are normal (see [Mat], Thm. 34). All this
yields that the relevant morphism is an isomorphism.
The lemma is proven.

\bigskip

As a corollary we get

\noindent{\bf Lemma 11}. Let $B_{F_0}$ be a smooth elementary cobordism.
The birational equivalence
$B_{F_0-}/K^*\mathrel{{-}\,{\rightarrow}} B_{F_0+}/K^*$
determined by a flip
\[\begin{array}{rccccc}
&B_{F_0-}/K^*&&&&B_{F_0+}/K^*\\
&&\psi_-\searrow&& \swarrow\psi_+&\\
&&&B_{F_0}//K^*&&\\
\end{array}\]
is the
composite of a toroidal blow-up
$\widetilde{B_{F_0-}/K^*\times_{B_{F_0}//K^*}B_{F_0+}/K^*
}\rightarrow B_{F_0-}/K^*$
and a toroidal blow-down
$B_{F_0+}/K^*\leftarrow
\widetilde{B_{F_0-}/K^*\times_{B_{F_0}//K^*} B_{F_0+}/K^*
}$ whose fibers are weighted projective spaces, where $\widetilde{B_{F_0-}/K^*\times_{B_{F_0}//K^*}B_{F_0+}/K^*
}$ is the
normalization of ${B_{F_0-}/K^*\times_{B_{F_0}//K^*}B_{F_0+}/K^*}$.

\noindent{\bf Proof}. This follows from Lemmas 9 and 10 and the analogous fact
for toric flips.

\bigskip

\noindent{\bf Theorem 1}. Let $\pi: X\rightarrow X'$ be a birational morphism of
smooth projective varieties defined over a field of
characteristic zero. Assume that $\pi $ is an isomorphism over
$U\subset X'$ .Then one can find a sequence  $X_0=X,X_1,..., X_k=X'$ of complete
varieties with cyclic singularities
together with morphisms $\pi_i:X_i\rightarrow X'$, which
are isomorphisms over $U$ such that for $i=0,...,k-1$ either

$\bullet X_{i+1}\rightarrow X_i $ is  a simple
toroidal blow-up   commuting with
$\pi_i$ and $\pi_{i+1}$  and whose fibers are
weighted projective spaces.
or

$\bullet X_{i+1}\leftarrow X_i $  is a simple toroidal blow-down   commuting with
$\pi_i$ and $\pi_{i+1}$ and whose fibers are
weighted projective spaces
or

$\bullet$ there exist a toroidal variety $Z_i$, a morphism $\pi_i^z:
Z_i\rightarrow X'$ and a diagram

 \[\begin{array}{rccccc}
&X_i&&&&X_{i+1}\\
&&\psi_i\searrow&& \swarrow\psi_{i+1}&\\
&&&Z_i&&\\
\end{array}\]

\noindent which is a simple
toroidal flip  commuting with
$\pi_i$, $\pi_i^z$ and $\pi_{i+1}$ respectively such that the fibers of
$\psi_i$ and
$\psi_{i+1}$  are
weighted projective spaces.

\bigskip
\noindent{\bf Proof.} By Proposition 2 we can find a smooth cobordism
$B=B(X,X')
$ over $X'$ which is  trivial over $U$.  Let
$F_0,..., F_k$ be connected fixed point set components such that
$F_i>F_j$ implies $i>j$. By
Proposition 1 and Lemma 9  we see that $B^{F_0...F_i}_-/K^*=
B^{F_0...F_{i-1}}_{F_i+}/K^*$
differs from $ B^{F_0...F_{i-1}}_-/K^*=B^{F_0...F_{i-1}}_{F_i-}/K^*$ for
$i=0,...,k$ by a
simple toroidal flip, a simple toroidal blow-up or a simple toroidal blow-down.

As a corollary  we get.

\bigskip

\noindent
{\bf Theorem
 2}. Let $X$ and $X'$ be smooth projective birationally
equivalent varieties defined over a field of
characteristic zero with isomorphic open subsets $U\subset X $
and $U'\subset X'$
. Then one can find a sequence of complete
varieties together with isomorphic open subsets
 $X_0=X\supset U_0=U,X_1\supset U_1,..., X_k=X'\supset U_k=U'$
 with cyclic singularities
such that for $i=0,...,k-1$ the birational equivalence $X_i\mathrel{{-}\,{\rightarrow}} X_{i+1}$
defines an isomorphism $U_i\simeq U_{i+1}$ and either

$\bullet X_{i+1}\rightarrow X_i  $ is  a simple
toroidal blow-up   whose fibers are
weighted projective spaces
or

$\bullet X_{i+1}\leftarrow X_i  $ is a simple toroidal blow-down  whose fibers are
weighted projective spaces.
or

$\bullet$ there exists a toroidal variety $Z_i$
 and a diagram

 \[\begin{array}{rccccc}
&X_i&&&&X_{i+1}\\
&&\psi_i\searrow&& \swarrow\psi_{i+1}&\\
&&&Z_i&&\\
\end{array}\]

 \noindent which is a simple
toroidal flip   such that the fibers of
$\psi_i$ and
$\psi_{i+1}$  are
weighted projective spaces.

\bigskip

\noindent {\bf Proof}. Let $\overline{X}$ be a smooth resolution
of the join $X*X'$. Apply Theorem 1 to the morphisms
$\overline{X}\rightarrow X$ and $\overline{X}\rightarrow X'$.

\bigskip

As a corollary  we get

\bigskip
\noindent {\bf Theorem 3}.
Let $X$ and $X'$ be smooth projective birationally
equivalent varieties defined over a field of
characteristic zero with isomorphic open subsets $U\subset X $
and $U'\subset X'$. Then one can find a sequence of complete toroidal quasismooth
varieties together with isomorphic open subsets
 $X_0=X\supset U_0=U,X_1\supset U_1,..., X_k=X'\supset U_k=U'$
 such that for $i=0,...,k-1$ the birational equivalence
$X_i\mathrel{{-}\,{\rightarrow}} X_{i+1}$ is either a toroidal blow-up or
toroidal blow-down
defining the isomorphism $U_i\simeq U_{i+1}$ whose fibers are
weighted projective spaces.

\noindent {\bf Proof}. Apply Lemma 11 to Theorem 2.

\bigskip
\noindent {\bf Remark}. By the Moishezon Theorem [Moi], which
says that each smooth complete variety over an algebraically closed field of
characteristic zero can be rendered projective by a sequence of
blow-ups with smooth centers, one can
prove   theorems  similar to Theorems 1, 2, 3 on smooth complete varieties.

\vspace*{1 cm}

\noindent {\bf REFERENCES}



\noindent {\bf [A]} Atiyah, M. F. {\it On an analytic surface
with double points}, Proc. Royal Soc. A 247 (1958), 237-244

\noindent {\bf [B-B]} Bia\l ynicki-Birula, A. {\it Some
theorems on actions of algebraic groups}, Annals of
Math.,98(1973), 480-497.

\noindent {\bf [B-B,S]} Bia\l ynicki-Birula, A. \'Swiecicka, J.
{\it Complete quotients by algebraic torus actions}. Lecture
Notes in Mathematics 956,
Springer-Verlag Berlin-Heidelberg-New York (1982)
,10-22.

\noindent {\bf [B,M]} Bierstone, E., Milman, P. {\it
Canonical desingularization in characteristic zero by
blowing-up of the maximum strata of a local invariant}.
Invent. Math. 128 (1997), no 2, 207-302.

\noindent {\bf [Dan]} Danilov,V.I. {\it The geometry of
toric varieties}, Russian Math. Surveys 33(1978),97-154;
Uspechi Mat. Nauk 33 (1978), 85-134

\noindent {\bf [D,H]} Dolgachev,I. Hu, Y. {\it Variation of
Geometric Invariant Theory Quotient}, Inst.Hautes \'Etudes Sci.
Publ.Math. No 87 (1998),5-56.

\noindent {\bf [Gro]} Grothendieck, A.
{\it EGA El\'ements de g\'eometrie alg\'ebrique}, Publ.Math. IHES.



\noindent {\bf [Har1]} Hartshorne, R. {\it Algebraic Geometry},
Springer-Verlag, New York, Heidelberg, Berlin, 1977

\noindent {\bf [Har2]} Hartshorne, R.
{\it Ample subvarieties of algebraic varieties}. Lecture Notes in
Mathematics 156,
Springer-Verlag Berlin-Heidelberg (1970).


\noindent {\bf [Hir]}
Hironaka H.{\it  Resolution of singularities of an
algebraic variety over a field of characteristic zero. I, II}
Annals of Math.79 (1964), 109-203, 205-326. MR 33:7333

\noindent {\bf [Jur]} Jurkiewicz J.{\it Torus embeddings,
polyhedra, $ k^*$-actions and homology}, \\ Dissertationes
Mathematicae CCXXXVI, 1985

\noindent {\bf [Kon]} Konarski J.{\it A pathological example of
$k^*$}, in: Group Action and Vector Fields, Lecture Notes in
Math.956, Springer-Verlag, Berlin 1982

\noindent {\bf [Lu]} Luna, D. {\it Slices \'etale}. Sur les
groupes alg\'ebriques, pp 81-105, Bull. Soc. Math. France, Paris,
Memoire 33 Soc.Math.France, Paris, 1973

\noindent {\bf [Mat]} Matsumura, H. {\it Commutative ring
theory.} Cambridge studies in Advanced Mathematics 8. Cambridge
University Press, Cambridge New York, 1988

\noindent {\bf [Mil]} Milnor, J.{\it Morse Theory}, Ann. of
Math. Stud. 51, Princeton Univ. Press, Princeton

\noindent {\bf [Moi]} Moishezon B.{\it On n-dimensional Compact
Varieties with n Algebraically Independent Meromorphic Functions},Amer
Math.Soc. Translation 63 (1967) 51-177.

\noindent {\bf [Mor]} Morelli, R. {\it The birational Geometry of
Toric Varieties}. J.Alg.Geom.5(1996) 751-782.

\noindent {\bf [Mum1]} Mumford D.{\it Geometric invariant theory}, Ergeb. Math.
Bd. 34, Springer Verlag, 1965.

\noindent {\bf [Mum2]} Mumford D.{\it Introduction To Algebraic Geometry}.

%
\noindent {\bf [Nag1]}
Nagata, M.
{\it Imbeddings of an abstract variety in a complete variety}. \\
 J.Math.Kyoto Univ.(1962), 1-10.

\noindent {\bf [Nag2]}
Nagata, M.
{\it On rational surfaces} I, Mem. Coll. Sci. Kyoto, A32 (1960), 351-370.

\noindent {\bf [Oda]} Oda, T. {\it Convex Bodies and Algebraic
Geometry}, Springer-Verlag, Berlin, Heidelberg, New-York,1988

\noindent {\bf [R]} Reid, M. {\it What is a flip?} Colloquium
Talk, Univ. of Utah, Dec 1992


\noindent {\bf [Sum]}
Sumihiro, H.
{\it Equivariant completion}. J.Math.Kyoto.Univ. 13(1974) 1-28.

\noindent {\bf [Tha1]} Thaddeus, M. {\it Stable pairs, linear
systems and the Verlinde formula.} Invent. Math. 117 (1994)
no.2, 317-353

\noindent {\bf [Tha2]} Thaddeus, M. {\it Toric quotients and
flips}. Topology, geometry and field theory, 193-213, World,
Sci., River Edge, NJ, 1994

\noindent {\bf [Tha3]} Thaddeus, M. {\it Geometric invariant
theory and flips}.
J. Amer.Math.Soc. 9 (1996) no. 3 , 691-723

\noindent {\bf [Wlo]} W\l odarczyk, J. {\it Decomposition of birational
toric maps in blow-ups and \\ blow-downs}. Trans.
Amer.Math.Soc.349(1997), no.1, 373-411.

\noindent {\bf [Zar]}
Zariski, O.
{\it Complete linear systems on normal varieties and a generalization
of a Lemma of Enriques-Severi}, Ann. of Math., 55 (1952) 552-592.


\end{document}